# A Hierarchical Multi-Objective Programming Approach to Planning Locations for Macro and Micro Fire Stations

Xinghan Gong[1], Jun Liang[1], Yiping Zeng[1], Fanyu Meng[1,2], Simon Fong[3], Lili Yang[1*]

[1]Department of Statistics and Data Science
Southern University of Science and Technology
Shenzhen, China

[2]Academy for Advanced Interdisciplinary Studies
Southern University of Science and Technology
Shenzhen, China

[3]Department of Computer and Information Science
University of Macau
Macau, China

e-mail: xg2356@columbia.edu, 12032857@mail.sustech.edu.cn, zengyp@sustech.edu.cn, mengfy@sustech.edu.cn, ccfong@um.edu.mo, yangll@sustech.edu.cn

*Abstract*—Fire stations are among the most crucial emergency facilities in urban emergency control system in terms of their quick response to fires and other emergencies. Location planning for fire stations has a significant influence on their effectiveness and capability of emergency responses trading off with the cost of constructions. To obtain efficient and practical siting plans for fire stations, various major requirements including effectiveness maximization, distance constraint and workload limitation are required to be considered in location models. This paper proposes a novel hierarchical optimization approach taking all the major requirements for location planning into consideration and bonds functional connections between different levels of fire stations at the same time. A single-objective and a multi-objective optimization model are established coupled with genetic algorithm (GA) with elitist reservation and Pareto-based multi-objective evolutionary algorithm for model solving. The proposed hierarchical location model is further performed in a case study of Futian District in Shenzhen, and the siting results justify the effectiveness and practicality of our novel approach.

## I. INTRODUCTION

The accelerated process of urbanization poses more and more challenges to the public emergency service in cities than ever before [1]. And there is an urgent need to optimally distribute emergency service facilities, including hospitals and fire stations, in response to the increase of large-scale emergencies [2]. Determination of locations to site emergency service facilities has a significant impact on people's lives and properties [3]. Fire stations, as primary fire rescue providers, are attached with great importance in city planning to ensure public safety and reduce property losses as well [4].

Typical fire station system in China consists of first class and second class fire stations, which are categorized into macro fire stations in this paper. Macro fire stations, as the main executive of fire rescue missions, hold primary responsibility to ensure fire safety in urban areas. However, the construction of macro fire stations consumes a large amount of land and financial resources in terms of setup and operation costs. Besides, macro fire stations are typically assigned with enormous service areas, resulting in a longer rescue time to the edge of service areas. Therefore, expanding only macro fire stations to enhance fire management system is no longer appropriate in current urban planning in regards of financial costs and time based rescue efficiency. The Chinese government proposed a hierarchical fire station structure in "Construction Standard of Urban Fire Stations 2017" ("Construction Standard" in the rest of the texts) [5] to fix this problem. In "Construction Standard", a novel type of fire station, micro fire station, is recommended to mainly serve local communities. Compared with macro fire stations, micro fire stations have a much lower setup and operation cost because of their less area occupancy and smaller service areas. The smaller service areas character further allows micro fire stations to have a much shorter rescue time, bringing high rescue efficiencies to local communities. Therefore, a reasonable location siting of macro fire stations coupled with micro fire stations is of great need for decision makers (DM) to strengthen the urban fire rescue forces.

Related papers show that multiple modeling methods have been conducted in urban fire station location problem, e.g. maximal covering problem (MCP), set covering problem (SCP) and p-median problem (PMP) [6]. Specifically, a fuzzy multi-objective programming combined with GA has been proposed to determine the optimal fire station locations to serve areas with different fire risk categories [7]. Yao *et al.* [8] constructed an bi-objective location model combined with SCP and PMP to locate fire stations adjacent to areas with high fire risks and minimize the total construction number as well. Aktaş [9] *et al.* applied both the SCP and MCP to acquire fire station location in 10 different scenarios in Istanbul. To make the optimal fire station locations more closely fit the actual situation, they employed Istanbul's road networks to determine additional fire station locations. Furthermore, a two ranked level model was developed to minimize response time to fire accidents with allocation strategies fully considered [10]. And Hendrick *et al.* also proposed that the influence of existing fire stations should not be neglected while considering where new ones are placed [11], [12].

*Corresponding author

Though differences exist in the location modeling techniques for fire station planning, some common principles and constraints are followed in these researches:

- Minimization of financial cost of fire stations. Minimizing the number of fire stations is often needed to meet the budget constraint or reduce the setup and operation costs [13], [14].
- Minimization of response time. Arrival assurance of fire engines at fire scenes within a certain time is the most basic standard the fire service authorities must follow [3], [7], [10].
- Maximization of total service areas. It is desired to cover as many demand points with high risks as possible [9].
- Minimization of the distance between fire stations and high-risk areas. With fire risks varying in space, optimal locations of fire stations should be close to high-risk regions to shorten the response time [15], [16].

Despite of the fact that hierarchical modeling method has been widely applied in emergency medical service systems and telecommunications networks problems [17], [18], [19], the existing researches in fire station planning mainly focused on modeling single leveled fire station system and few studies paid much attention to the current hierarchical structure of macro and micro fire stations. Yu *et al.* [3] adopted MCP and SCP models to arrange the locations of macro and micro fire stations with GA applied into solving the hierarchical model. A basic assumption lay in [3] is to relocate all macro fire stations instead of adding new ones. However, in practical application, it is unrealistic and of high costs to site macro fire stations in this way. Besides, insufficient factors have been considered in [3], as two single objective optimization problems with limited constraints are proposed for macro and micro fire station locations. Therefore, this paper aims to construct a more realistic and practical hierarchical location model for macro and micro fire stations based on the government's requirements in "Construction Standard". The mathematical representations for macro and micro fire stations models are single objective programming and multi-objective programming, respectively. Our hierarchical also involves the effect of existing fire stations and real road network in potential location points for fire stations. A case study in Futian District, Shenzhen indicates that our modeling method has strong portability and is applicable to the emergency management authorities.

The main contributions of this paper include:

- This paper provides a practical construction approach for cities without a complete multi-level fire facility system.
- This paper combines the covering location model with geographic information system (GIS) approach to better simulate the real-world situation. And the visualization makes our modelling procedures more understandable.
- Considering the hierarchical location structure of macro and micro fire stations, we take the following aspects into our model construction:
  – Take factors including transportation and rescue efficiency into account.
  – Apply multi-objective programming in the second-level model.
  – Establish an interaction between fire stations with different levels.

The remainder of this paper is organized as follows. Section II presents our hierarchical location model for macro and micro fire stations. Section III conducts a case study in Futian District, Shenzhen, China, with GA methods applied to solve the model. The results validate the practicality and efficiency of our solving and modeling approaches. Then, the paper is ended with some conclusions in section IV.

## II. HIERARCHICAL LOCATION MODEL

In this section, a hierarchical model is proposed to solve the location problem of macro and micro fire stations. In specific, an extended MCP with distance constraint between neighboring stations, and a combination model of SCP and PMP with capacity and distance constraints are adopted to select locations for macro and micro fire stations, respectively. The following subsections introduces the two models in detail.

### A. Model Parameters

The following parameters will be applied into our hierarchical location model.

- $I$: Set of demand point, $i \in I$;
- $J_1$: Set of potential macro fire station locations, $j_1 \in J_1$;
- $J_2$: Set of potential micro fire station locations, $j_2 \in J_2$;
- $d_{ij_k}$: Distance between $i$ and $j_k$, $k=1,2$;
- $d_{jj'}$: Distance between two adjacent fire stations, $j, j' \in J_k$, $k=1,2$. For any two stations located at $j_a, j_b \in J_k, k=1,2$, we call the fire station located at $j_a$ "adjacent" to the another fire station located at $j_b$ if the distance between these two stations is shorter than the distance between the fire station at $j_a$ and any other fire stations [7]. The detailed definition of adjacent fire stations is described as follows:

$$d_{j_a j_b} \leq d_{jj_a}, \forall j \in J_k, j \neq j_a, j_b, k=1,2. \quad (1)$$

- $d^{s_1}, d^{l_1}$: Minimal and maximal distance between each pair of adjacent macro fire stations, respectively;

- $d^{s_2}, d^{l_2}$ : Minimal and maximal distance between each micro fire station and its adjacent macro fire station, respectively;
- $a_i$ : Demand value of node $i$. In this study, demand value is defined as the estimated fire risk value;
- $N$ : Number of macro fire stations to be located;
- $S$ : Workload limitation value of micro fire station. In this paper, the workload is defined as the total sum of fire risks covered by each fire station;
- $\Omega_{ij_k} = \{j_k \mid d_{ij_k} \leq R_k, k=1,2\}$: Set of potential locations $j_k$ capable of serving demand point $i$;
- $\eta_{ij_k} = \{i \mid d_{ij_k} \leq R_k, k=1,2\}$: Set of demand points $i$ that can be covered by potential fire station at $j_k$;
- $\Phi$: Set of existing macro fire station locations;
- $R_1, R_2$ : Service radius of macro and micro fire stations, respectively;
- 
$$X_{j_k} = \begin{cases} 1, & \text{if a potential fire station is sited at } j_k \\ 0, & \text{else} \end{cases}$$

- 
$$Y_{ij_k} = \begin{cases} 1, & \text{if node } i \text{ is suitably covered by station at } j_k \\ 0, & \text{else} \end{cases}$$

### B. Macro Fire Station Model

In our macro fire station model, we propose to keep all the existing macro fire stations and build a limited number of new macro stations to cover as many high-risk communities as possible. And the distance constraints between adjacent stations are adopted to obtain a balanced distribution of newly-built macro stations.

*1) Extended MCP for Macro Fire Stations:*

The proposed model is formulated as follows,
- Objectives:

$$\text{Max} \sum_{i \in I} a_i Y_{ij_1} \quad (2)$$

- Subject to :

$$\sum_{j_1 \in \Omega_{ij_1}} X_{j_1} \geq Y_{ij_1}, \quad \forall i \in I, \quad (3)$$

$$\sum_{j_1 \in J_1} X_{j_1} = N, \quad (4)$$

$$d^{s_1} \leq d_{j_1 j_1'} \leq d^{l_1}, \quad \forall j_1, j_1' \in J_1, \quad (5)$$

$$d^{s_1} \leq d_{j_1 j_1'} \leq d^{l_1}, \quad \forall j_1 \in J_1, j_1' \in \Phi, \quad (6)$$

$$X_{j_1}, Y_{ij_1} = \{0,1\}, \quad \forall i \in I, j_1 \in J_1. \quad (7)$$

Equation (2) is an objective aiming to maximize the total demand coverage of communities by macro fire stations, thereby encouraging macro fire station to cover as many high -risk communities as possible. Constraint (3) ensures that each chosen community would be covered by at least one macro fire station. Constraint (4) specifies the number of new macro fire stations to locate. Constraint (5) limits the distance between each pair of new adjacent macro fire stations. Furthermore, if a new macro fire station is adjacent to an existing macro fire station, constraint (6) limits the distance between these two fire stations. Finally, constraint (7) imposes binary coding requirement on decision variables $X_{j_1}$ and $Y_{ij_1}$.

### C. Micro Fire Station Model

Micro fire stations have a low setup and operation cost and are convenient to relocate because of their assembled building structure. Therefore, even if the initial allocation of micro fire stations is unreasonable, DM still can relocate micro fire stations with low costs.

In current fire management system, micro fire station is taken as a supplementary rescue force for macro fire stations. In specific, since the macro fire station model is to cover as many high-risk communities as possible, there exists some communities with smaller fire risks that cannot be covered by macro fire stations. Therefore, micro fire stations are sited to cover these communities and provide auxiliary support to macro fire stations in high-risk communities as well. Furthermore, service area of each micro fire station is required to intersect with at least one macro fire station's service area. Therefore, based on these siting requirements, we combine SCP and PMP with a novel objective in maximizing the average distance of adjacent micro fire stations to construct our model.

*1) Spatial Location Model for Micro Fire Stations:* The proposed model is formulated as follows,
- Objectives:

$$\text{Min} \sum_{j_2 \in J_2} X_{j_2} \quad (8)$$

$$\text{Min} \sum_{i \in I} \sum_{j_2 \in \Omega_{ij_2}} a_i d_{ij_2} Y_{ij_2} \quad (9)$$

$$\text{Max} \frac{\sum_{j_2, j_2' \in J_2} d_{j_2 j_2'}}{|J_2|} \quad (10)$$

- Subject to:

$$\sum_{j_2 \in \Omega_{ij_2}} X_{j_2} \geq 1, \quad \forall i \in I, \quad (11)$$

$$Y_{ij_2} \leq X_{j_2}, \quad \forall i \in I, j_2 \in J_2, \quad (12)$$

$$\sum_{i \in \eta_{j_2}} a_i \leq S, \quad \forall i \in I, \quad (13)$$

$$d^{s2} \leq d_{j_1 j_2} \leq d^{l2}, \quad \forall j_1 \in J_1, j_2 \in J_2, \quad (14)$$

$$X_{j_2}, Y_{ij_2} = \{0,1\}, \quad \forall i \in I, j_2 \in J_2. \quad (15)$$

The objective (8) is to minimize the total number of newly-built micro fire stations. The objective (9) is to minimize the distance between micro fire stations and communities with higher fire risks. And objective (10) is to maximize the average distance of each pair of adjacent micro fire stations. Constraint (11) and (12) ensure that all the communities can be covered by at least one micro fire stations. Constraint (13) limits the workload of each micro fire station where we define the total sum of fire risk values within the service area as the workload value. Constraint (14) limits the distance between each micro station and its adjacent macro fire station. And finally, binary requirements are imposed in constraint (15).

The combination of the objectives and constraints constitute a multi-objective decision optimization model for the micro fire station location problem, shown as

$$\text{Min } [F1, F2, F3], \quad (16)$$

where $F1 = Obj1$ (8), $F2 = Obj2$ (9) and $F3 = -Obj3$ (10).

## III. CASE STUDY

The proposed model is implemented and applied in an empirical study evaluating fire service in Futian District, Shenzhen, China. The interest is to construct a new fire station system while considering the connection between different levels of fire stations and the existing fire stations' effects. Our research scenario here is to keep all the existing macro fire stations and determine optimized locations for new macro and micro fire stations. The model is solved by using Python Geatpy library [20], and all the experiments are run on a computer with an Intel Core i7 2.59 GHz CPU and 8.0 GB RAM. We also combine commercial GIS software, ArcGIS (version 10.6) into selecting potential fire station location points and visualizing our results.

### A. Study Area

The research area, Futian District, is located in Shenzhen, Guangdong Province, China. Shenzhen is one of the mega-cities in southern China and Futian District is the central area of Shenzhen, which takes up, according to "Shenzhen Statistical Yearbook-2020" [21], 16.88% (454.65 billion Yuan) of Shenzhen total GDP (2692.71 billion Yuan) and has a population density of 20769 person per square kilometers. With such a high economic status and density of population, Futian District is undertaking great pressure from fire control management. Therefore, there exists an urgent need for the city's DM to enhance the fire rescue forces in the region.

### B. Data Preparation

In the case study, we obtain our potential location points based on the road network constructed from the open street map (OSM). And in our calculation of estimated fire risk value of each community, we acquire the historical fire accident location points between 2014 and 2019 in Futian District from the Emergency Management Bureau of Shenzhen Municipality, and the population density value from the open platform of Shenzhen government [21]. In total, there are 1789 historical fire accident points, 4 existing macro fire station and 95 communities in Futian District. And our novel hierarchical location model is utilized to generate siting locations for macro and micro fire stations.

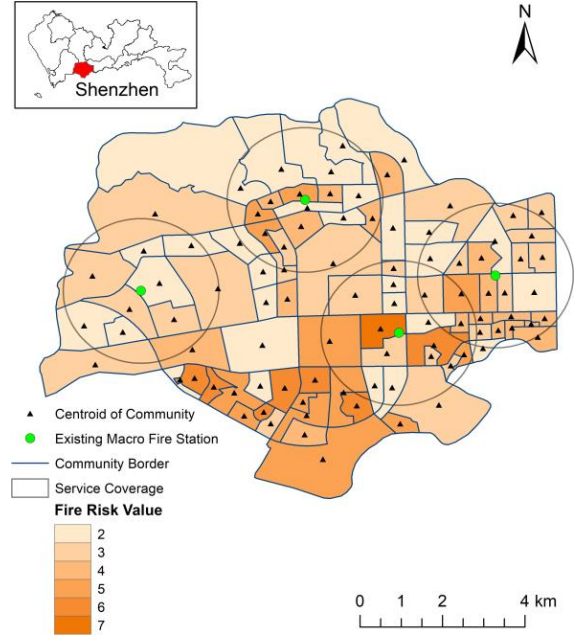

Fig. 1. Study area, existing macro fire stations and fire risk values.

### C. Location Siting of Macro Fire Stations

In "Construction Standard", the government requires that the response time of macro fires station should be within 4 minutes, and the scale of fire stations' service area is also limited for the service area of first-class fire station should not exceed 7 $km^2$ and second-class fire station's under 4 $km^2$. Besides, "Construction Standard" provides the follow equation to calculate service area of fire station

$$A = 2P^2, \quad (17)$$

where $A$ is the service area of fire station and $P$ is the coverage radius of service area. Thus, taken these regulations and real road structure into account, we propose the following weighted method to get the radius of service area in our model

$$R_1 = \sqrt{\frac{1}{2}(\beta \cdot A_1 + (1-\beta) \cdot A_2)}, \; \beta \in (0,1), \quad (18)$$

where $A_1$ and $A_2$ are the recommended service areas of first and second class fire stations respectively based on "Construction Standard". $R_1$ represents the service radius of macro fire station in the model, and $\beta$ is the weight value determined by DM. In the case study, to reflect the priority of first class fire station and combined with fire service authorities' opinions, the value of $\beta$ is set to 0.7 and with

$A_1 = 7\ km^2$ and $A_2 = 4\ km^2$. Then, the service radius of macro fire station $R_1$ in the model is set as 1.746 km. In this section, we introduce the key steps to obtain the locations of new macro fire stations as follows:

*1) Calculation of Estimated Fire Risk Values:* In the hierarchical model, the demand value $a_i$ is defined as the estimated fire risk value of each community. By this way, in the first level of the hierarchical model, newly-built macro fire stations are required to locate near high-risk communities while satisfying all the constraints. Thus, we first need to calculate the estimated fire risk value.

TABLE I
THE FIRE RISK VALUE OF COMMUNITIES IN TERMS OF THE NUMBER OF HISTORICAL FIRE ACCIDENTS

| The number of historical fire accidents | The number of communities | Risk value $r_a$ |
|---|---|---|
| $2 < n \leq 14$ | 45 | 1 |
| $14 < n \leq 26$ | 33 | 2 |
| $26 < n \leq 46$ | 10 | 3 |
| $46 < n \leq 99$ | 7 | 4 |

TABLE II
THE FIRE RISK VALUE OF COMMUNITIES IN TERMS OF POPULATION DENSITY

| Population density ρ | The number of communities | Risk value $r_p$ |
|---|---|---|
| $0.26 < \rho \leq 3.60$ | 53 | 1 |
| $3.60 < \rho \leq 8.02$ | 28 | 2 |
| $8.02 < \rho \leq 13.37$ | 12 | 3 |
| $13.37 < \rho \leq 24.85$ | 2 | 4 |

Regional fire risk in urban communities is closely interrelated with population size and the frequency of fire accidents within local communities. Therefore, in the case study, we combine the population density with the number of fire accidents in each community to calculate the estimated fire risk value. First, we apply natural breaks [22] in our model to rank fire risk levels of communities. Natural breaks has advantage in searching natural breakpoints which have statistical meanings in sequences, and with these obtained breakpoints, the sequence can be classified into groups with similar properties. As presented in Table I and II, we rank the fire risk levels of communities from population density and total number of fire accidents with each aspect dividing into four categories using natural breaks method. Then, we combine them with a weight parameter $\gamma$, $\gamma \in (0,1)$, to get the final estimated fire risk values of communities according to equation

$$a_i = \gamma \cdot r_a + (1-\gamma) \cdot r_p. \quad (19)$$

For general purposes, we adopt the parameter $\gamma$ as 0.5 [3]. we then map fire risk values into each corresponding community in Fig. 1, where communities with deeper colors have higher fire risk values. Specifically, Fig. 1 shows the study area, estimated fire risk and the location of existing macro fire stations. In the case study, centroid points are taken from main residential area of each community and we assume that fire stations can serve communities as long as their service areas can cover their centroids. It can be seen from Fig. 1 that existing macro fire stations have a unbalanced distribution for communities with deeper colors in the southern area of Futian District cannot be served by these fire stations. Therefore, there exists an urgent need to construct new macro fire stations in this region.

*2) Calculation of the Number of New Macro Fire Stations:* Tzeng *et al.* [16] and Yang *et al.* [7] calculate the number of new fire stations to be built by balancing the total setup and operating cost of fire stations and total lost cost of incidents in a given area. The detailed expression of the optimal cost model is summarized as follows,

$$\text{Min } f(N) = N \times SC + \alpha \times TLC \times e^{-N}, \quad (20)$$

where $SC$ is the setup and operating cost of each fire station, $TLC$ is the total lost cost of incidents in a given area and $N$ is the number of new fire stations to be built. And by setting the derivative of $f(N)$ as zero, the number of new fire stations $N$ is given as

$$N = \text{int}(\log TLC - \log SC + \beta). \quad (21)$$

However, for the initial lost cost model in [16] with $\alpha = 1$ in (20), the value of $SC$ is much higher than $TLC$ in megacities (e.g. Shenzhen) during normal years based on the data from open data platform of Shenzhen government [21]. As a result, this causes $N < 0$ in (21) and cannot be applied into our model. As for the lost cost model with $\alpha > 1$, the historical record utilized to calculate $\alpha$ is normally hard to obtain because of the loss of relative data in fire management system. Therefore, our model proposes to obtain the number of new macro fire stations to be built based on the coverage area of macro fire station and total area of given research region with

$$N = \left\lceil \frac{TAR - N_e \times SAM}{SAM} \right\rceil, \quad (22)$$

where $TAR$ is the total area of research region, $N_e$ is the number of existing macro fire stations and $SAM$ is the service area of each macro station. In (22), we round up the value because the equation does not consider in the overlapping of service area in terms of the irregular-shape character of research regions which would lead to the unavailability of the actual number of new macro fire stations. Hence, the actual number of new macro fire stations is higher than the result of (22) and we round up the value to solve this problem.

In the case study, the total area of Futian District is 78.66 $km^2$ and urban area takes up 75% of it (i.e. 58.995 $km^2$). Also,

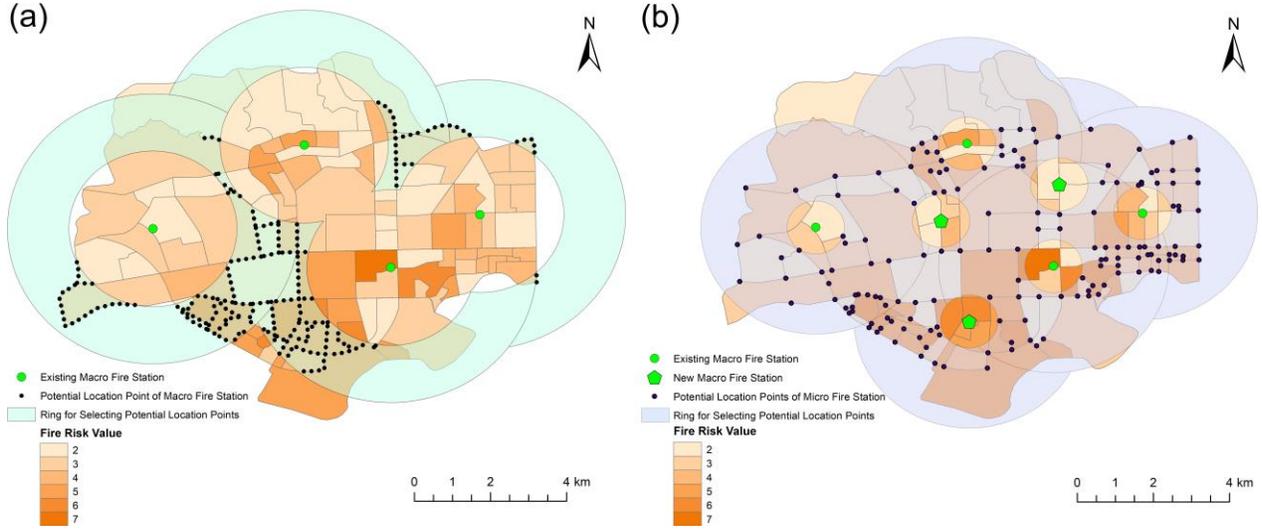

Fig. 2. Final potential location points for macro fire stations (a) and micro fire stations (b).

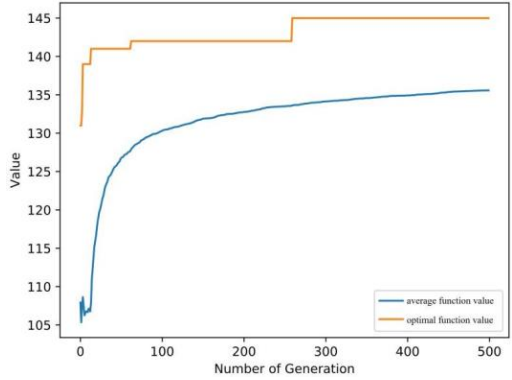

Fig. 3. Total sum of fire risk values of elite individual in each generation.

with the service area limitation in "Construction Standard" and our model assumptions, the service area of each macro fire station is 9.58 $km^2$. Then, according to (22), 3 new macro fire stations need to be sited in the model.

*3) Selection of Potential Location Points for Macro Fire Stations:* The division of communities in Futian District is based on road network, and according to "Construction Standard", locations of fire stations are required to have high transportation accessibility. Therefore, in the case study, we first find all intersections in the road network of Futian District and take these intersection points as potential location points for macro fire stations. Specifically, since the number of newly-built macro fire stations has been settled in section III-C2, we propose to expand choices of potential location points for macro fire stations by taking additional points every 200 meters within each road segment. However, since each road segment is not exact multiple of 200, we delete potential location points within 50 meters of road intersections which have higher priorities than other points because of their transportation convenience. Besides, as shown in Fig. 2(a), in order to fulfill the distance limitation between each new macro fire station and its adjacent existing macro fire station, we take final choices of potential location points within the ring centered at location points of existing macro fire stations with a width of $d^{l_1} - d^{s_1}$, where

$$d^{l_1} = 2R_1 + \varepsilon, \tag{23}$$

and $d^{s_1}$ is obtained by requiring the overlapping area of adjacent macro fire stations cannot exceed 30% of the total service area. In (23), each pair of adjacent macro fire stations should not be too far from each other with a tolerance of $\varepsilon = 0.05\ km$. And we obtain total 239 potential locations points for new macro fire stations in Fig. 2(a).

*4) Application of GA:* Finally, we apply GA with elitist reservation [20] to solve the single objective extended MCP model of macro fire stations in (2)-(7). And as shown in Fig. 3, the total sum of fire risk value gradually increases and convergent in the 270th generation. The final siting result is presented in Fig. 4(a), and if we classify community with fire risk value equal or larger than 4 as high-risk community, the newly-built macro fire stations solve the problem of unbalanced distribution of existing macro fire stations in southern region of Futian District where gathered with high-risk communities and cover more communities as well. Besides, the location siting result of macro fire stations can in total cover 87.18% (i.e. 34 communities) of all 39 high-risk communities and can also cover 88.43% of all the 1789 historical fire accidents in Futian District. In conclusion, with constructing three new macro fire stations, we can arrange the distribution of macro fire stations more efficient and enhance fire rescue forces of the fire management system as well.

### D. Location Siting of Micro Fire Stations

According to "Construction Standard", micro fire stations target at serving communities and with a limitation of 3

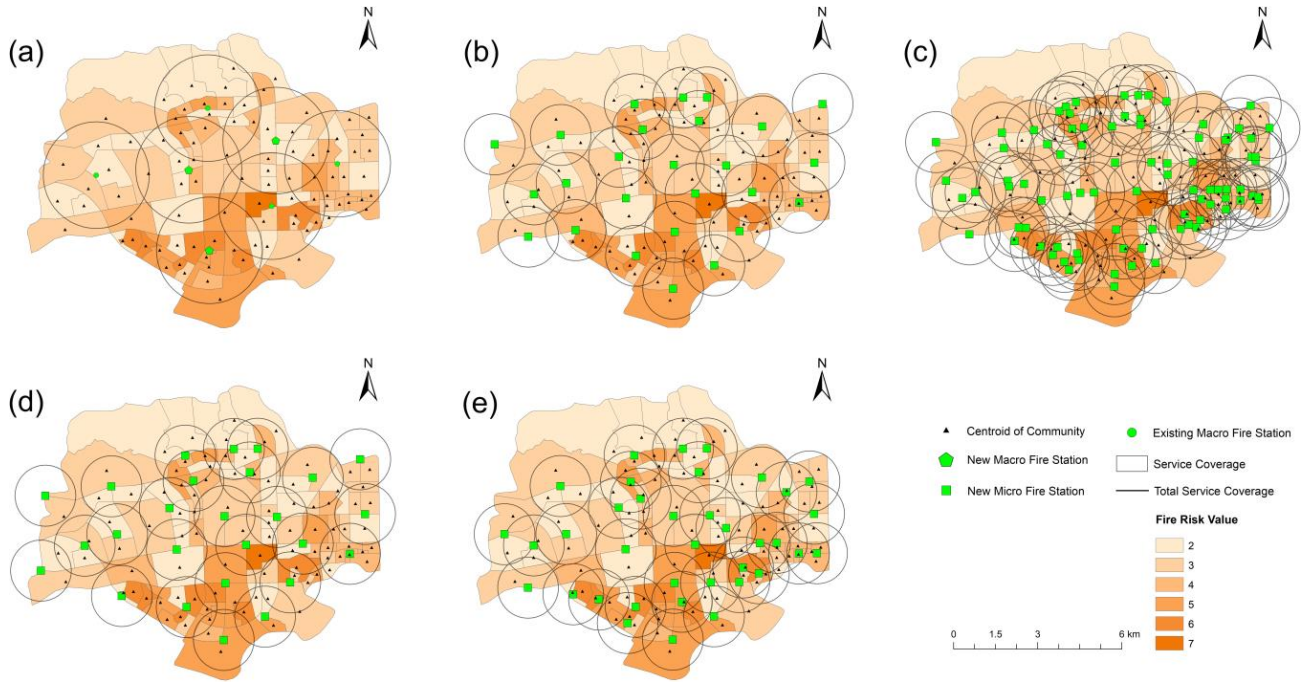

Fig. 4. Siting results of macro fire stations (a) and four representative siting solutions of micro fire stations: A is preferred by F1 (b); B is preferred by F2(c); C is preferred by F3 (d); and D is the equally weighted solution (e).

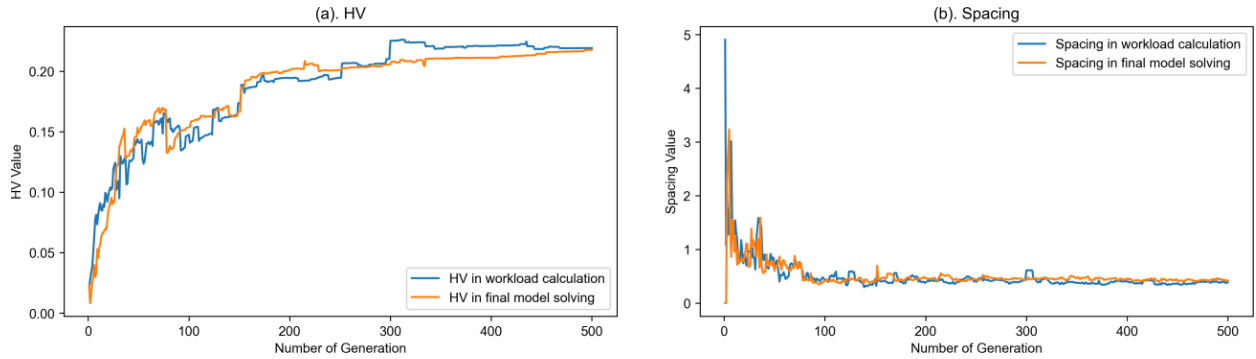

Fig. 5. The values of HV (a) and Spacing (b) in section III-D2 and III-D3.

minutes response time and service area under 1 $km^2$. Then, based on (17), service radius of micro fire station is set to 1 $km$. The demand value $a_i$ of communities is the same with the fire risk value obtained in section III-C1. The key steps on obtaining siting results of micro fire stations are introduced as follows:

*1) Selection of Potential Location Points:* In the case study, the potential location points of micro fire stations are also obtained from the road network in Futian District. The difference is that we only take intersection points of road network this time. This is because the number of micro fire stations to be sited is not fixed and with too many potential location points, convergence process of the GA applied to solve the model will be too time-consuming, leading to a low efficiency to get the siting results. So, similar method in section III-C3, we obtain final potential location points within the ring which centers at all the macro fire stations including the existing ones and the results got from the macro fire station model, and with a width of $d^{l_2} - d^{s_2}$, as presented in Fig. 2(b). The values of $d^{l_2}$ and $d^{s_2}$ are obtained from (24) and (25) with a tolerance value $\varepsilon = 0.05\ km$ respectively,

$$d^{l_2} = R_1 + R_2 + \varepsilon, \qquad (24)$$
$$d^{s_2} = R_1 - R_2 - \varepsilon. \qquad (25)$$

Through this step, we can ensure the service area of each micro fire station is adjacent to a macro fire station's service area, and there are total 165 final potential location points selected for micro fire stations in Fig. 2(b).

*2) Calculation of Workload Limitation Value:* Among various models in solving the location problem of fire stations, a usual assumption is that facilities are without workload limitations which means as long as the location of a fire accident is sited within the service of a fire station, it will be suitably resolved by the fire station. This assumption might work with macro fire stations for their high construction standard with more rescue equipment. However, micro fire stations do not possess such high rescue capability and a micro fire station with excessive communities to manage will lead to unbalancedness and inefficiency in fire management system [3]. Therefore, in the second-level of our model, we add a capacity constraint (13) so as to balance the number of high-risk communities can be covered by micro fire stations. And the workload is working as the upper average total fire risk value which can be covered by each micro fire station. The calculation of workload value is shown as follows,

$$S = \frac{\sum_{i=1}^{M} \max_i \{S_{ij}\}}{M}, \ 1 \leq j \leq |D_i|, \ i, j \in \mathbb{Z}, \quad (26)$$

where $S_{ij}$ is the $j$ th workload value from $i$ th run. $D_i$ is the siting solution space of $i$ th run and $M$ is the number of runs.

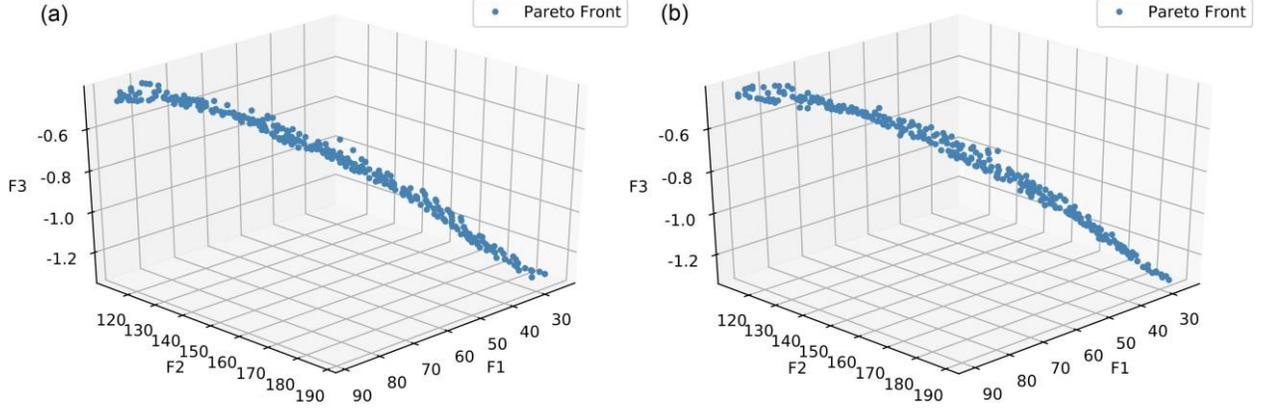

Fig. 6. Best final Pareto optimal fronts among the 10 runs: in workload calculation section (a); in final model solving section (b).

To obtain the value of $S_{ij}$, we first take away the workload constraint (13) in the model and acquire the siting results preferred by F2. This is because F2 encourages each micro fire station to cover as many high-risk communities as possible and thereby making each micro fire station reach its upper average total fire risk coverage value. Then $S_{ij}$ is the average sum of fire risk values of each siting result preferred by F2.

In the case study, the model without constraint (13) is solved by using NSGA-II. NSGA-II is a highly efficient MOEA and has been widely applied into solving multi-objective optimization problems [23]. NSGA-II involves an elitist-preserving approach to speed up the performance of the algorithm and help prevent the loss of elitist solutions once they are found. Furthermore, a fast non-dominated sorting approach is applied to reduce the computation complexity, and a crowding distance evaluation mechanism is used to preserve the diversification of Pareto-optimal solutions. Detailed framework of NSGA-II can refer to [24].

Encoding scheme in MOEA also has great effect on the performance of the algorithm [25]. In the case study, unlike [3], [23] and [26] which adopt a mesh-based spatial representation strategy, we obtain our potential location points on the road segments. Therefore, we apply the classical binary encoding scheme with "1" representing the locations to be sited with micro fire stations and "0" representing the otherwise in NSGA-II.

Due to the stochastic nature of evolutionary algorithms, it is necessary to perform several parallel runs to evaluate their performance. Therefore, without loss of generality, 10 independent runs (i.e. $M = 10$) are implemented by using fixed parameter values but with different initial populations [23]. In specific, the maximum number of evolutionary generations is 500. The size of the initial population is limited to 300, the crossover probability is 0.9 and the mutation probability is 0.005. The best final Pareto-optimal front among 10 runs is displayed in Fig. 6(a). As presented in Fig. 6(a), the result shows the ability of NSGA-II in finding a well-convergent Pareto-optimal front.

The performance of NSGA-II is evaluated based on hypervolume (HV) and Spacing indicators. The HV indicator is described as the volume of the space in the objective space dominated by the Pareto front approximation and delimited from above by a reference point [27]. It has been utilized to evaluate the convergence and diversity of MOEA and the bigger the HV is, the batter the performance is. In the case study, we take the nadir point [28] in each generation as the reference point and detailed calculation algorithm of HV can refer to [29]. The Spacing indicator captures the variation of the distance between elements of a Pareto front approximation and a lower value of the MOEA is considered to be better [27]. This paper adopts the Spacing metric to evaluate the distribution character of non-dominated solution points along each Pareto front approximation.

First, we analyze the convergence performance of NSGA-II based on HV metric. The values of HV for every generation are shown in Fig. 5(a). We find that the HV values sharply increase for the first 100 generations and finally become stable for 200 generations from 300th to 500th generation, which shows that NSGA-II is able to converge around 300th generation. Furthermore, as seen in Fig. 5(b), non-dominated solution points are able to uniformly distributed along Pareto front approximations from around 100th generation because the spacing values keep stable for 400 generations from 100th generation. It shows the capability of NSGA-II in generating a set of uniformly distributed Pareto optimal solutions. Therefore, it is appropriate to apply NSGA-II into solving our model and the final computation result of workload value $S$ according to (26) is 19.374829.

*3) Model Solution and Analysis:* The spatial location model for micro fire stations with capacity and distance constraints (8)-(15) is also solved by NSGA-II. Similarly, 10 independent runs are implemented with the same parameters defined in section III-D2 and the best result is shown in Fig. 6(b). As presented in Fig. 6(b), it shows the ability of NSGA-II in finding a well-convergent Pareto-optimal front approximation.

The model verification is also based on HV and Spacing metrics. For the evaluation of convergence performance of the algorithm, as shown in Fig. 5(a), NSGA-II converges from around 200th generation for the HV values keep stable for 300 generations from 200th to 500th generation. As for the distribution property of solutions, Fig. 5(b) represents that the Pareto solution sets are able to uniformly distributed starting from 100th generation.

TABLE III
OBJECTIVE VALUES AND INCIDENT COVERAGE RATES OF THE FOUR SITING SOLUTIONS

| Solution | A | B | C | D |
|---|---|---|---|---|
| F1 | 26 | 86 | 26 | 37 |
| F2 | 187.87 | 113.98 | 188.91 | 140.75 |
| F3 | -1.31 | -0.42 | -1.33 | -0.85 |
| Incident Coverage Rate | 93.63% | 99.27% | 94.47% | 96.81% |

The 275 non-dominated solutions on the final Pareto optimal front of NSGA-II together constitute a candidate pool for the DM, where siting solutions of micro fire stations with different preferences can be chosen as the final location plans. However, there is no best solution among these non-dominated solutions because they cannot dominate each other by definition. In the case study, an equally weighted solution and three extremely preferred solutions are selected among the 275 non-dominated solutions [26]. The four kinds of solutions are respectively labeled as A, B, C and D. For each solution, statistical analyses are presented in Table III and the result visualization is displayed in Fig. 4.

Fig. 4(b) represents solution "A" which is preferred by F1. The locations of the 26 micro fire stations have an overall uniform distribution with more concentration around the northern high-risk communities of Futian District. However, not all high-risk communities, especially communities in the southern areas, receive attention from the location solution and each micro fire station located in the southern areas obviously has to cover excessive high-risk communities, which shows a unbalanced distribution of the solution. This is because solution "A" only taking the number of stations to be constructed as the most important factor. Fig. 4(c) represents solution "B" which is preferred by F2. The number of selected locations is 86 and the solution achieved the distribution of more concentration around all high-risk communities. However, the high overlap rate around these high-risk communities shows an economic inefficiency of the solution and therefore, the solution "B" is not recommended to the DM. Fig. 4(d) represents solution "C" which is preferred by F3. As shown in Fig. 4(d), solution "C" has the most uniform distribution because of the maximal distance between each pair of adjacent micro fire stations. Similar defect in solution "A", southern micro fire stations is overloaded with high-risk communities to cover, resulting in a unbalanced distribution. Equally weighted solution "D" shown in Fig. 4(e) has the most balanced trade-off among financial limitation, rescue capability and rescue efficiency among all the non-dominated solutions. As presented in Fig. 4(e), the number of selected locations is 37 and all the high-risk communities are concentrated with more micro fire stations. Therefore, solution "D" is recommended to the DM as the final siting location plan.

In conclusion, 3 new macro fire stations and 37 micro fire stations obtained from the whole hierarchical model together construct the new fire station system in Futian District. There are three major advantages in our result. First, a connection between different levels of fire stations is constructed in the model, making our results more reasonable in reality. Second, with keeping all the existing macro fire stations and constructing a limited number of new ones, the siting results allows macro fire stations to cover more high-risk communities and micro fire stations, as a supplementary rescue forces in fire management system, cover the remaining ones. Finally, the results fulfill our initial aim to require micro fire stations to reach the fire accident scene at the very beginning while waiting for the adjacent macro fire stations to implement major rescue mission. More specifically, the coverage rate of high-risk communities is 100% and the entire fire station system could also cover 98.66% (i.e.1765) of total historical fire accidents in Futian District as well. Therefore, with our hierarchical location model, a reasonable and efficient location siting instruction is able to be provided to DM.

IV. CONCLUSIONS

Considering the multi-level structure of fire station system proposed by the Chinese government, this paper constructs a hierarchical model for macro and micro fire stations. In specific, in the first level of the hierarchical model, an extended MCP with distance constraint between adjacent stations is selected for macro fire stations siting because of the fixed number of new macro fire stations to be sited. Then, in the second level of our hierarchical model, a combination model of SCP and PMP with capacity and distance constraints is adopted for micro fire stations, as micro fire stations are established to mainly serve local communities. The hierarchical model is further performed in a case study of Futian District in Shenzhen, and a non-dominated sorting genetic algorithm II (NSGA-II) is proposed to solve the model. NSGA-II is a highly efficient and elitist-preserving multi-objective evolutionary algorithm (MOEA). In the case study, we adopt hypervolume (HV) and Spacing indicators to evaluate the performance of NSGA-II, and the test results show its feasibility in solving the hierarchical model. Furthermore, the final siting results of macro and micro fire stations prove the effectiveness of our model, which can be utilized to generate valid reference for decision makers.

ACKNOWLEDGMENT

This research was supported by National Key R&D Program of China (2019YFC0810700 and 2018YFC-0807000), National Natural Science Foundation of China (71771113, 71704091 and 71804026 No.72004141) and Shenzhen Science and Technology Plan Project(N0. JSGG20180717170802038) and Basic and Applied Basic Research Foundation of Guangdong Province (No.2019-A1515111074).

# Authors' background

1. This form helps us to understand your paper better, the form itself will not be published.
2. Title can be chosen from: master student, Phd candidate, assistant professor, lecture, senior lecture, associate professor, full professor

| Your Name | Position | Research Field | Personal Webpage |
|---|---|---|---|
| Xinghan Gong | master student | Multi-objective optimization, Emergency management. | |
| Jun Liang | master student | Operations research, Risk modeling. | |
| Yiping Zeng | assistant professor | Public safety and emergency management, Natural language process, Pedestrian and evacuation dynamics. | https://stat-ds.sustech.edu.cn/teacher/ZENGYiping |
| Fanyu Meng | assistant professor | Risk analysis for natural disaster and accident, Emergency management and decision making, Road safety, Economic impact analysis for public health events, Driving psychology and behaviors. | https://faculty.sustech.edu.cn/mengfy/en/ |
| Simon Fong | associate professor | E-Commerce, Data mining, Business intelligence, Intelligent agent technology, Electronic governance. | https://www.fst.um.edu.mo/people/ccfong/#contact |
| Lili Yang | full professor | Big data analysis, Decision science, Multi-objective optimization, Emergency model building, Energy relationship analysis. | https://faculty.sustech.edu.cn/yangll/en/ |